\documentclass[a4paper,11pt,twoside]{article}
\usepackage{pst-fill,pst-grad}
\usepackage{textcomp}
\usepackage[english]{babel}
\usepackage[utf8x]{inputenc}
\usepackage{graphicx}
\usepackage{amsmath}
\usepackage{multicol}
\usepackage{float}
\usepackage{fancyhdr}
\usepackage[matrix,arrow,curve]{xy}
\usepackage{pstricks} 
\usepackage{amsmath,amsfonts,verbatim,afterpage,theorem,euscript,mathrsfs,amssymb}
\usepackage{amsfonts}
\usepackage{amssymb}
\usepackage{array}
\usepackage{dsfont}
\usepackage{hyperref}
\usepackage{color}

\def \vu{\vec{u}}

\def \vU{\vec{U}}
\def \vB{\vec{B}}
\def \vA{\vec{A}} 
\def \vUU{\vec{\mathcal{U}}}
\def \vBB{\vec{\mathcal{B}}}
\def \vV{\vec{V}}
\def \vW{\vec{W}}

\def \vn{\vec{\nabla}}
\def \vb{\vec{b}}
\def \vf{\vec{f}}
\def \vg{\vec{g}}
\def \vomega{\vec{\omega}}
\def \vrho{\vec{\rho}}
\def \vphi{\vec{\varphi}}
\def \vpphi{\vec{\phi}}

\def \R{\mathbb{R}^{3}}
\def \ds{\displaystyle}

\newtheorem{Proposition}{Proposition}[section]
\newtheorem{Lemme}{Lemma}[section]
\newtheorem{Theoreme}{Theorem}
\newtheorem{Corollaire}{Corollary}[section]
\newtheorem{Remarque}{Remark}[section]

\title{\bf  On the local regularity theory for the MHD equations}

\usepackage{authblk}
\author[1]{D. Chamorro\footnote{Corresponding author: \emph{diego.chamorro@univ-evry.fr}}}
\author[2]{F. Cortez}
\author[3]{J. He}
\author[4]{O. Jarr\'in}
\affil[1,3]{\footnotesize LaMME, Universit\'e d'Evry Val d'Essonne, France.}
\affil[2]{\footnotesize Escuela Polit\'ecnica Nacional, Quito, Ecuador.}
\affil[4]{\footnotesize Universidad T\'ecnica de Ambato, Ambato, Ecuador.}
\begin{document}
\maketitle
\begin{scriptsize}
\abstract{Local regularity results are obtained for the MHD equations using as global framework the setting of parabolic Morrey spaces. Indeed, by assuming some local boundedness assumptions (in the sense of parabolic Morrey spaces) for weak solutions of the MHD equations it is possible to obtain a gain of regularity for such solutions in the general setting of the Serrin regularity theory. This is the first step of a wider program that aims to study both local and partial regularity theories for the MHD equations.}\\
		
\end{scriptsize}
\section{Introduction}
In this article we study \emph{local} regularity results for the incompressible 3D magnetohydrodynamic (MHD) equations which are given by the following system:
\begin{equation}\label{EcuacionMHD}
\begin{cases}
\partial_{t}\vu=\Delta \vu -(\vu\cdot\vn)\vu+(\vb\cdot\vn)\vb-\vn p+\vf,\quad div(\vu)=0,\\[3mm]
\partial_{t}\vb=\Delta \vb -(\vu\cdot\vn)\vb+(\vb\cdot\vn)\vu+\vg,\quad div(\vb)=0,\\[3mm]
\vu(0,x)=\vu_{0}(x) \mbox{ \; and \; } \vb(0,x)=\vb_{0}(x), \,\,div(\vu_0)=0,\,\, div(\vb_0)=0, 
\end{cases}
\end{equation}
where $\vu, \vb:[0,T]\times \R\longrightarrow \R$ are two divergence-free vector fields which represent the velocity and the magnetic field, respectively, and the scalar function $p:[0,T]\times \R\longrightarrow \mathbb{R}$ stands for the pressure. The initial data $\vu_0, \vb_0: \R\longrightarrow \R$ and the external forces $\vf, \vg:[0,T]\times \R\longrightarrow \R$ are given and for the external forces we will always assume that they belong to the space $L^2_tH^1_x$.\\

Now if $\Omega\subset [0, +\infty[\times \R$ is a bounded set, we will say that the couple $(\vu, \vb)\in L^\infty_tL^2_x\cap L^2_tH^1_x(\Omega)$ satisfy the MHD equations (\ref{EcuacionMHD}) in the weak sense if for all $\vphi, \vpphi \in \mathcal{D}(\Omega)$ such that $div(\vphi)=div(\vpphi)=0$, we have\\
$$
\begin{cases}
\langle\partial_{t}\vu-\Delta \vu +(\vu\cdot\vn)\vu-(\vb\cdot\vn)\vb-\vf|\vphi \rangle_{\mathcal{D}'\times \mathcal{D}}=0,\\[3mm]
\langle\partial_{t}\vb-\Delta \vb +(\vu\cdot\vn)\vb-(\vb\cdot\vn)\vu-\vg|\vpphi \rangle_{\mathcal{D}'\times \mathcal{D}}=0,
\end{cases}
$$
note that if $(\vu, \vb)$ are solutions of the previous system, then there exists a pressure $p$ such that (\ref{EcuacionMHD}) is fulfilled in $\mathcal{D}'$.\\

It is clear that if the magnetic field $\vb=0$, then the previous equations (\ref{EcuacionMHD}) are reduced to the classical Navier-Stokes equations 
\begin{equation}\label{NSequations}
\partial_{t}\vu=\Delta \vu -(\vu\cdot\vn)\vu-\vn p+\vf, \qquad div(\vu)=0, 
\end{equation}
for which some results related to \emph{regularity} are available. Indeed, let us briefly recall the Serrin regularity theory for the classical Navier-Stokes system: 
\begin{Theoreme}[local regularity, \cite{Serrin1}]\label{Theo_Serrin}
Let $Q=]a,b[\times B(x_{0},r_{0})$ be a bounded set where $]a,b[$ is an interval and $B(x_{0},r_0)$ is a ball with $x_{0}\in \mathbb{R}^{3}$ and $r_0>0$. Let $\vf\in L^2_{t}H^k_{x}(Q)$ for some $k\geq 0$, let $\vu\in L^{\infty}_{t}L^{2}_{x}(Q)\cap L^{2}_{t}\dot{H}^{1}_{x}(Q)$ and $p\in \mathcal{D}'(Q)$; if we assume that $\vu$ is a weak solution on $Q$ of the Navier-Stokes equations (\ref{NSequations}) then, if 
\begin{equation}\label{HypotheseSerrin}
\vu\in L^{\infty}_{t}L^{\infty}_{x}(Q),
\end{equation}
we obtain that locally the regularity of $\vu$ is given by the regularity of the external force $\vf$:  for every $a<c<b$ and $0<\rho<r_0$ we have that $\vu \in L^{\infty}\big(]c,b[, H^{k+1}(B(x_{0},\rho))\big)\cap L^{2}\big(]c,b[,\dot{H}^{k+2}(B(x_{0},\rho))\big)$. The points of $]0, +\infty[\times\R$ for which we have the condition (\ref{HypotheseSerrin}) for some $Q$ will be called \emph{regular points}.
\end{Theoreme}
Remark that no particular assumption is needed for the pressure $p$, which can be a very general object and this fact is a very important feature of this theory. 
\begin{Remarque}\label{RemarqueSerrinLpLq}
Note that the assumption $\vu\in L^{\infty}_{t}L^{\infty}_{x}(Q)$ stated in (\ref{HypotheseSerrin}) can be generalized. Serrin \cite{Serrin1} proved that, if $\vf \in L^2_{t}H^1_{x}(Q)$ and if 
\begin{equation}\label{HypoSerrinLpLq}
\vu\in L^p_{t} L^q_{x}(Q)\quad \mbox{with } \tfrac{2}{p}+\tfrac{3}{q}<1, 
\end{equation}
then for every $a<c<b$ and $0<\rho<r_{0}$ we have that $\vu\in L^\infty_{t} L^\infty_{x}(]c,b[\times B(x_{0}, \rho))$.  
\end{Remarque}
Important and significant efforts have been made to generalize even more this hypothesis (\ref{HypotheseSerrin}), see for example \cite{Struwe}, \cite{Taka} or \cite{Chen}. In particular, parabolic Morrey-Campanato spaces were used by O'Leary \cite{OLeary}, see also \cite{PGLR2}, to generalize Serrin's theorem and we will see how to exploit this framework for the MHD equations (\ref{EcuacionMHD}).\\

It is worth to mention here that another regularity theory is available for the Navier-Stokes equations. Indeed, Caffarelli, Kohn and Nirenberg developed in \cite{CKN} a second approach, known as the \emph{partial regularity theory}, which is essentially based on energy estimates. Of course these two points of view (local and partial) are quite different since they require different hypotheses\footnote{In particular, for the Caffarelli-Kohn-Nirenberg theory some information is needed for the pressure $p$, which is not the case for the Serrin theory.} and since the results obtained are obviously different, however -and this point is important- a common treatment can be performed by using the framework of parabolic Morrey spaces. See for example Kukavica \cite{Kukavica} for generalization of the Caffarelli-Kohn-Nirenberg theory in this parabolic setting. One special feature of this common framework appears to be crucial when studying the role of the pressure in the Caffarelli-Kohn-Nirenberg theory for the classical Navier-Stokes equations, indeed, as it is shown in \cite{CML}, the language of parabolic Morrey spaces is a powerful tool which allows to mix, in a very specific sense, these two regularity theories.\\

Although many studies concerning regularity are available for the MHD equations (see for example \cite{Chen}, \cite{Jia} or \cite{Larios} and the references there in for a generalization of Theorem \ref{Theo_Serrin} and Remark \ref{RemarqueSerrinLpLq} to the MHD equations), to be best of our knowledge, a detailed treatment using parabolic Morrey spaces is missing. Since this framework is important to improve the understanding of the role of the pressure in these regularity theories, we find interesting to set up in this article the first step of our approach -given by Theorem \ref{TheoPrincipal} below- that will eventually lead us in a forthcoming work to define new classes of solutions for the MHD equations (\ref{EcuacionMHD}).\\

The plan of the article is the following. In Section \ref{Secc_NotationAndPresentation} we introduce some notation and we present our main theorem while in Section \ref{Secc_PropertiesMorrey} we recall some useful fact about parabolic Morrey spaces. Finally, in Section \ref{Secc_ProofTh} we detail the proof of all the results stated before. Some classical but useful results are gathered in the appendix. 
\section{Notation and presentation of the results}\label{Secc_NotationAndPresentation}
Before stating the main theorem of this article, we need to introduce some notation related to parabolic Morrey spaces. It is worth noting here that the use of these parabolic spaces is actually given by the underlying structure of the MHD equations: indeed, in one hand we have that if $(\vu, p, \vb)$ is a solution of (\ref{EcuacionMHD}), then for $\lambda>0$ the triplet $\lambda \vu(\lambda^2t, \lambda x)$, $\lambda^2 p(\lambda^2t, \lambda x)$ and $\lambda \vb(\lambda^2t, \lambda x)$ is still a solution of the MHD equations and this remark will lead us to a very particular dilation structure. On the other hand, when studying existence for these equations, we can see the system (\ref{EcuacionMHD}) as a nonlinear perturbation of the heat equation and thus the properties of the heat kernel $h(\sqrt{t},x)$ must also to be taken into account. It is thus natural to consider the homogeneous space $(\mathbb{R}\times \R, d, \mu)$ where $d$ is the parabolic quasi-distance given by 
\begin{equation}\label{Def_QuasiDistance}
d\big((t,x), (s,y)\big)=|t-s|^{\frac{1}{2}}+|x-y|,
\end{equation}
and where $\mu$ is the usual Lebesgue measure $d\mu=dtdx$. Remark that the homogeneous dimension is now $Q=5$. See \cite{Folland} for more details concerning the general theory of homogeneous spaces. \\

Now for $1< p\leq q<+\infty$, parabolic Morrey spaces $M_{t,x}^{p,q}$ are  defined as the set of measurable functions $\vphi:\mathbb{R}\times\R\longrightarrow \R$ that belong to the space $(L^p_tL^p_x)_{loc}$ such that $\|\vphi\|_{M_{t,x}^{p,q}}<+\infty$ where
\begin{equation}\label{DefMorreyparabolico}
\|\vphi\|_{M_{t,x}^{p,q}}=\underset{x_{0}\in \R, t_{0}\in \mathbb{R}, r>0}{\sup}\left(\frac{1}{r^{5(1-\frac{p}{q})}}\int_{|t-t_{0}|<r^{2}}\int_{B(x_{0},r)}|\vphi(t,x)|^{p}dxdt\right)^{\frac{1}{p}}.
\end{equation}
Remark that we have $M_{t,x}^{p,p}=L^p_{t}L^p_x$. In Section \ref{Secc_PropertiesMorrey} we will present some useful properties of these spaces. \\

As we are interested in studying \emph{local} regularity properties of the solutions of the MHD equations (\ref{EcuacionMHD}), in what follows we will always consider here the following subset of $]0,+\infty[\times\R$:
\begin{equation}\label{DefConjuntoOmega}
\Omega=]a,b[\times B(x_{0},r), \quad \mbox{with} \quad 0<a<b<+\infty, x_{0}\in \R \mbox{ and } 0<r<+\infty. 
\end{equation}
The main theorem of this article reads as follows.
\begin{Theoreme}\label{TheoPrincipal}
Let $\vu_{0}, \vb_{0}:\R\longrightarrow \R$ such that $\vu_{0}, \vb_{0}\in L^{2}(\R)$ and $div(\vu_{0})=div(\vb_{0})=0$ be two initial data and consider two external forces $\vf, \vg:[0,+\infty[\times\R\longrightarrow \R$ such that $\vf, \vg\in L^{2}([0,+\infty[, \dot{H}^{1}(\R))$. \\
Assume that $p\in \mathcal{D}'(\Omega)$ and that $\vu, \vb:[0,+\infty[\times\R\longrightarrow \R$ are two vector fields that belong to the space
\begin{equation}\label{UBSolutionFaiblesMHD}
L^{\infty}(]a,b[, L^{2}(B(x_{0},r)))\cap L^{2}(]a,b[, \dot{H}^{1}(B(x_{0},r))),
\end{equation}
such that they satisfy the MHD equations (\ref{EcuacionMHD}) over the set $\Omega$ given in (\ref{DefConjuntoOmega}).\\
		
\noindent If moreover we have the following local hypotheses
\begin{equation}\label{LocalHypo1}
\begin{cases}
\mathds{1}_{\Omega}\vu\in M_{t,x}^{p_{0},q_{0}}(\mathbb{R}\times\R) & \mbox{ with } 2<p_{0}\leq q_{0}, 5<q_{0}<+\infty\\[3mm]
\mathds{1}_{\Omega}\vb\in M_{t,x}^{p_{1},q_{1}}(\mathbb{R}\times\R) &\mbox{ with } 2<p_{1}\leq q_{1}, 5<q_{1}<+\infty, 
 \end{cases}
\end{equation} 
and $p_1 \leq p_0$, $q_1 \leq q_0$, then, for all $\alpha, \beta$ such that $a<\alpha<\beta<b$ and for all $\rho$ such that $0<\rho<r$, we have 
$$\vu\in L^{q_0}(]\alpha,\beta[, L^{q_0}(B(x_{0},\rho)))\quad \mbox{and}\quad \vb\in L^{q_1}(]\alpha,\beta[, L^{q_1}(B(x_{0},\rho))).$$
\end{Theoreme}
Note that once we have this result -and observing that the parameters above satisfy the condition (\ref{HypoSerrinLpLq})- we can thus apply Remark \ref{RemarqueSerrinLpLq}  to obtain that we actually have $\vu, \vb\in (L_t^{\infty}L_x^\infty)_{loc}$ and then, by the Serrin theory stated in Theorem \ref{Theo_Serrin} in the context of the MHD equations \cite{ChenMiao}, we will deduce local regularity for the solutions of the MHD equations.
\section{Useful properties of parabolic Morrey spaces}\label{Secc_PropertiesMorrey}
We state here some results that will be frequently used in the sequel. The first one is just a consequence of H\"older's inequality.
\begin{Lemme}\label{Lemme_Product}
If $\vf, \vg:\mathbb{R} \times \R\longrightarrow \R$ are two function that belong to the space $M_{t,x}^{p,q} (\mathbb{R} \times \R)$ then we have the inequality
$$\|\vf\cdot\vg\|_{M_{t,x}^{\frac{p}{2}, \frac{q}{2}}}\leq  C\|\vf\|_{M_{t,x}^{p, q}} \|\vg\|_{M_{t,x}^{p, q}}.$$
\end{Lemme} 
Our next lemma explains the behaviour of parabolic Morrey spaces with respect to localization in time and space. 
\begin{Lemme}\label{Lemme_locindi}
Let $\Omega$ be a bounded set of $\mathbb{R} \times \R$ of the form given in (\ref{DefConjuntoOmega}). If we have $1< p_0 \leq p_1$, $1< p_0\leq q_0 \leq q_1<+\infty$ and if the function $\vf:\mathbb{R} \times \R\longrightarrow \R$  belongs to the space $M_{t,x}^{p_1,q_1} (\mathbb{R} \times \R)$ then we have the following localization property 
$$\|\mathds{1}_{\Omega}\vf\|_{M_{t,x}^{p_0, q_0}} \leq C\|\vf\|_{M_{t,x}^{p_1,q_1}}.$$
\end{Lemme} 
Let us now introduce, for $0<\alpha < 5$, the parabolic Riesz potential $\mathcal{I}_{\alpha}$ of a locally integrable function $\vf:\mathbb{R} \times \R\longrightarrow \mathbb{R}^3$  which is given by the expression
\begin{equation}\label{Def_Riesz_potential}
\mathcal{I}_{\alpha} (\vf) (t,x) =\int_{\mathbb{R}} \int_{\R} \frac{1}{(\vert t-s\vert^{\frac{1}{2}} + \vert x-y \vert)^{5-\alpha}} \vf (s,y) dy\, ds.
\end{equation}
As for the standard Riesz Potential in $\mathbb{R}^3$, we have a corresponding  boundedness property:
\begin{Lemme}[Adams-Hedberg's inequality for parabolic Riesz potentials]\label{Lemme_Hed}
If $0 <\alpha < \frac{5}{q}$, $1 < p \leq q < +\infty$ and $\vf \in M_{t,x}^{p,q} (\mathbb{R} \times \R)$ then for $\lambda=1- \frac{\alpha q}{5}$ (which verifies $0<\lambda<1$) we have the inequality
$$\|\mathcal{I}_\alpha(\vf)\|_{M_{t,x}^{\frac{p}{\lambda} , \frac{q}{\lambda}}} \leq C\|\vf\|_{M_{t,x}^{p,q}}.$$
\end{Lemme}
See \cite{Adams} for a proof of this fact. From these general lemmas we will now deduce two specific results that will be helpful in our computations. 
\begin{Corollaire}\label{Coro_I1}
Let $\Omega$ be a bounded set of the form given in (\ref{DefConjuntoOmega}). If $2 < p \leq q$, $5 < q \leq 6$, and $\vf \in M_{t,x}^{\frac{p}{2}, \frac{q}{2}}(\mathbb{R} \times \R)$, then we have
\begin{itemize}
\item[1)] $ \mathds{1}_{\Omega} \mathcal{I}_1(\vf) \in M_{t,x}^{\frac{p}{\lambda}, \frac{q}{\lambda}}(\mathbb{R} \times \R),$ with $\lambda = 1 - \frac{q-5}{5q}$ (remark that $0<\lambda<1$).
\item[2)] $ \mathds{1}_{\Omega} \mathcal{I}_1(\vf) \in M_{t,x}^{\delta, q}(\mathbb{R} \times \R), $ where $\delta = \text{min} (\frac{p}{\lambda}, q)$ with the same $\lambda$ as before.
\end{itemize}
\end{Corollaire}
{\bf Proof.} For the first point, it is enough to notice that since $2 < p \leq q$ and $5 < q \leq 6$ we have $2\sigma\leq \lambda$ where $\lambda = 1 - \frac{q-5}{5q}$ and $\sigma=1-\frac{q}{10}$. Thus, applying Lemmas \ref{Lemme_locindi} and \ref{Lemme_Hed} we have:
$$\|\mathds{1}_{\Omega}\mathcal{I}_1(\vf)\|_{M_{t,x}^{\frac{p}{\lambda}, \frac{q}{\lambda}}}\leq C\|\mathcal{I}_1(\vf)\|_{M_{t,x}^{\frac{p}{2\sigma}, \frac{q}{2\sigma}}}\leq C \|\vf\|_{M_{t,x}^{\frac{p}{2}, \frac{q}{2}}}.$$
For the second point, since we have $\delta = \text{min} (\frac{p}{\lambda}, q) \leq \frac{p}{\lambda}$ and $q < \frac{q}{\lambda}$, by Lemma \ref{Lemme_locindi} we can write $\|\mathds{1}_{\Omega}\mathcal{I}_1(\vf)\|_{M_{t,x}^{\delta, q}}\leq \|\mathds{1}_{\Omega}\mathcal{I}_1(\vf)\|_{M_{t,x}^{\frac{p}{\lambda}, \frac{q}{\lambda}}}$ and it only remains to apply the first point just proved.\hfill$\blacksquare$
\begin{Corollaire}\label{Coro_I2}
Let $\Omega$ be a bounded set of the form given in (\ref{DefConjuntoOmega}). If $2 < p \leq q$, $5 < q \leq 6$ and $\vf \in M_{t,x}^{\frac{p}{2}, \frac{q}{2}} (\mathbb{R} \times \R)$, then we have 
$$ \mathds{1}_{\Omega} \mathcal{I}_2 (\mathds{1}_{\Omega} \vf) \in M_{t,x}^{\delta, q} (\mathbb{R} \times \R), $$
where $\delta = \text{min} (\frac{p}{\lambda}, q)$ with $\lambda = 1 - \frac{q-5}{5q}$.
\end{Corollaire}
{\bf Proof.} Notice first that we cannot use Lemma \ref{Lemme_Hed} directly since we are dealing here with the Riesz potencial $\mathcal{I}_{\alpha}$ with $\alpha=2 > \frac{5}{q/2}$. To overcome this gap we will exploit the double localization of the function $\mathds{1}_{\Omega} \mathcal{I}_2 (\mathds{1}_{\Omega} \vf)$. Indeed, observing that $\delta=\text{min} (\frac{p}{\lambda}, q)\leq q$, we can write by Lemma \ref{Lemme_locindi}
$$\|\mathds{1}_{\Omega} \mathcal{I}_2 (\mathds{1}_{\Omega} \vf) \|_{M_{t,x}^{\delta, q}}\leq C\|\mathds{1}_{\Omega} \mathcal{I}_2 (\mathds{1}_{\Omega} \vf) \|_{M_{t,x}^{q, q}}.$$
Consider now a parameter $\sigma$ such that $\sigma < \frac{5}{2} < \frac{q}{2}$ and such that $\sigma$ is close enough to $\frac{5}{2}$ so that we have
$\frac{\sigma}{1- 2\sigma/5} \geq \frac{\text{min} (\frac{p}{2}, \sigma)}{1- 2\sigma/5} > q$, thus we have by Lemma \ref{Lemme_locindi} $\|\mathds{1}_{\Omega} \mathcal{I}_2 (\mathds{1}_{\Omega} \vf) \|_{M_{t,x}^{q, q}}\leq C\| \mathcal{I}_2 (\mathds{1}_{\Omega} \vf) \|_{M_{t,x}^{\frac{\text{min} (\frac{p}{2}, \sigma)}{1- 2\sigma/5}, \frac{\sigma}{1- 2\sigma/5}}}$. Since now we do have the condition $2< \frac{5}{\sigma}$, by Lemma \ref{Lemme_Hed} we deduce the inequality
$$\| \mathcal{I}_2 (\mathds{1}_{\Omega} \vf) \|_{M_{t,x}^{\frac{\text{min} (\frac{p}{2}, \sigma)}{1- 2\sigma/5}, \frac{\sigma}{1- 2\sigma/5}}}\leq C\| \mathds{1}_{\Omega} \vf \|_{M_{t,x}^{\text{min} (\frac{p}{2}, \sigma), \sigma}}.$$
It is enough to remark that $\text{min} (\frac{p}{2}, \sigma)\leq \frac{p}{2}$ and that $\sigma \leq \frac{q}{2}$ to obtain $\| \mathds{1}_{\Omega} \vf \|_{M_{t,x}^{\text{min} (\frac{p}{2}, \sigma), \sigma}}\leq C\| \vf \|_{M_{t,x}^{\frac{p}{2}, \frac{q}{2}}}$ and the Corollary \ref{Coro_I2} follows. \hfill$\blacksquare$

\section{Proof of Theorem \ref{TheoPrincipal}} \label{Secc_ProofTh}
The first thing to do is to define our framework, thus from a general parabolic ball $\Omega$ of the type (\ref{DefConjuntoOmega}) that will be fixed once and for all, we consider the two following subsets: 
\begin{equation}\label{Def_ParabolicCylinders}
\Omega_0= ]\alpha, \beta[\times B(x_{0}, \rho) \quad \mbox{and}\quad \Omega_1= \left]\frac{a+\alpha}{2}, \frac{b+\beta}{2}\right[\times B\left(x_{0}, \frac{r+\rho}{2}\right),
\end{equation}
and remark that since $0<a<\alpha<\beta<b$ and $0<\rho<r$, we have the inclusion
\begin{equation}\label{InclusionsOmegas}
\Omega_0 \subset \Omega_1 \subset \Omega.
\end{equation}
Note in particular that the conclusion of Theorem \ref{TheoPrincipal} is given over the subset $\Omega_0$. \\

Observe also that since we are working in a local setting, due to the localization property stated in Lemma \ref{Lemme_locindi} and with no loss of generality we may assume in hypothesis (\ref{LocalHypo1}) that we have $5< q_0, q_1 < 6$.\\

Once our framework is clear, in order to prove Theorem \ref{TheoPrincipal} we will use the following strategy: we define two technical parameters $0<\lambda_{0}, \lambda_{1}<1$ such that
\begin{equation}\label{Def_Lambda01}
\lambda_{0}=1-\frac{q_{0}-5}{5q_{0}}\quad \mbox{and}\quad \lambda_{1}=1-\frac{q_{1}-5}{5q_{1}},
\end{equation}
and we prove that we have
\begin{equation}\label{maingoal_ub}
\mathds{1}_{\Omega_0}\vu\in M_{t,x}^{\sigma_{0},q_{0}}(\mathbb{R}\times\R), \quad
\mathds{1}_{\Omega_0}\vb\in M_{t,x}^{\sigma_{1},q_{1}}(\mathbb{R}\times\R),
\end{equation}
where $\sigma_0=\min\{\frac{p_{0}}{\lambda_{0}},q_{0}\}$ and $\sigma_1=\min\{\frac{p_{1}}{\lambda_{1}},q_{1}\}$. Now if  \eqref{maingoal_ub} holds, then by iteration we will be able to obtain 
\begin{equation}\label{EstimateIteration1}
\mathds{1}_{\Omega_0}\vu\in M_{t,x}^{q_{0},q_{0}} = L^{q_0}_t L^{q_0}_x, \quad \mathds{1}_{\Omega_0}\vb\in M_{t,x}^{q_{1},q_{1}} = L^{q_1}_t L^{q_1}_x,
\end{equation}
which is the conclusion of Theorem \ref{TheoPrincipal}. Indeed, if we have at our disposal estimates (\ref{maingoal_ub}), then reapplying the same arguments we will obtain $\mathds{1}_{\Omega_0}\vu\in M_{t,x}^{\sigma_{(0,1)},q_{0}}(\mathbb{R}\times\R), \quad
\mathds{1}_{\Omega_0}\vb\in M_{t,x}^{\sigma_{(1,1)},q_{1}}(\mathbb{R}\times\R)$, where $\sigma_{(0,1)}=\min\{\frac{\sigma_{0}}{\lambda_{0}},q_{0}\}=\min\{\frac{p_{0}}{\lambda_{0}^2},q_{0}\}$ and $\sigma_{(1,1)}=\min\{\frac{\sigma_{1}}{\lambda_{1}},q_{1}\}=\min\{\frac{p_{1}}{\lambda_{1}^2},q_{1}\}$, then observing that we have $\lim\limits_{n\to +\infty} \frac{p_{0}}{\lambda^n_{0}} = +\infty$ and $\lim\limits_{n\to +\infty} \frac{p_{1}}{\lambda^n_{1}} = +\infty$, we obtain  (\ref{EstimateIteration1}).\\

Now, to prove \eqref{maingoal_ub} we introduce two test functions $\phi,\varphi:\mathbb{R}\times \R\longrightarrow \mathbb{R}$ that belong to the space $\mathcal{C}^{\infty}_{0}(\mathbb{R}\times \R)$ and such that 
\begin{eqnarray}
\phi\equiv 1 \; \text{on} \; \Omega_0\quad \text{and} \quad \text{supp}(\phi)\subset \Omega_1,\label{DefSoporteFuncTest1}\\[2mm]
\varphi\equiv 1 \; \text{on}\; \Omega_1 \quad\text{and} \quad \text{supp}(\varphi) \subset  \Omega.\label{DefSoporteFuncTest2}
\end{eqnarray}
These functions satisfy two important facts: first we have $\phi(0,\cdot)=\varphi (0, \cdot) = 0$ and second due to the inclusions stated in (\ref{InclusionsOmegas}) we have the identity $\phi \varphi\equiv\phi$ in the whole space.\\
	
We define now $\vU=\phi\vu$ and $\vB=\phi\vb$. As long as we are interested in the behavior of $\vu$ and $\vb$ inside the set $\Omega_0$ and with the properties of the localization functions $\phi$ and $\varphi$ defined above, we can write
\begin{eqnarray*}
\vU&=&\varphi\left(\frac{1}{\Delta}\Delta (\phi \vu)\right)=\varphi\left(\frac{1}{\Delta}\left(\phi\Delta \vu-(\Delta\phi)\vu+2\sum_{i=1}^{3}\partial_{i}\big((\partial_{i}\phi)\vu\big)\right)\right),\\
\vB&=&\varphi\left(\frac{1}{\Delta}\Delta (\phi \vb)\right)=\varphi\left(\frac{1}{\Delta}\left(\phi\Delta \vb-(\Delta\phi)\vb+2\sum_{i=1}^{3}\partial_{i}\big((\partial_{i}\phi)\vb\big)\right)\right).
\end{eqnarray*}
Thus, verifying (\ref{maingoal_ub}) amounts to prove that $\vU\in M_{t,x}^{\sigma_{0},q_{0}}$ and $\vB\in M_{t,x}^{\sigma_{1},q_{1}}$ and for this we will first study in the expressions above the terms where the Laplacian does not act directly over the functions $\vu$ and $\vb$. More precisely if we define the quantities
\begin{equation}\label{Def_VW}
\vV=\varphi\left(\frac{1}{\Delta}(\phi\Delta \vu)\right)\quad \mbox{and}\quad \vW=\varphi\left(\frac{1}{\Delta}(\phi\Delta \vb)\right),
\end{equation}
we will study in the next lemma the behavior of the quantities $\vU-\vV$ and $\vB-\vW$ and we will prove that locally they belong to the parabolic Morrey spaces we are looking for.
\begin{Proposition} \label{Propo_UmoinsV}
Under the notation (\ref{Def_Lambda01}), assume $5<q_0<6$ and let $\sigma_0=\min\{\frac{p_{0}}{\lambda_{0}},q_{0}\}$, then we have $\mathds{1}_{\Omega}(\vU-\vV) \in  M_{t,x}^{\sigma_0,q_0}(\mathbb{R}\times\R)$. Symmetrically, if $5<q_1<6$ and if $\sigma_1=\min\{\frac{p_{1}}{\lambda_{1}},q_{1}\}$ then we have  $\mathds{1}_{\Omega}(\vB-\vW) \in M_{t,x}^{\sigma_1,q_1}(\mathbb{R\label{key}}\times\R)$.
\end{Proposition}
{\bf Proof of Proposition \ref{Propo_UmoinsV}.} 
We claim first that 
\begin{equation}\label{claim}
\vU-\vV= \varphi\left( \frac{1}{\Delta} \left( -(\Delta \phi) \vu + 2\sum_{i=1}^{3}\partial_{i}\big((\partial_{i}\phi)\vu\big) \right) \right)\in  L^{\infty}(]0, +\infty[, L^{6}(\R)).
\end{equation}
Indeed, recall that $\vu \in L^{\infty}(]a,b[, L^{2}(B(x_{0},r)))$ hence  $\vu \in L^{\infty}(]a,b[, L^{\frac{6}{5}}(B(x_{0},r)))$ and by definition of the test function $\phi$ we have $(\Delta \phi) \vu \in L^{\infty}(]0, +\infty[, L^{\frac{6}{5}}(\R))$ thus, recalling that we have by duality the embedding $L^{\frac{6}{5}}\subset \dot{H}^{-1}$, we obtain
$$(\Delta \phi) \vu \in L^{\infty}(]0, +\infty[, \dot{H}^{-1}(\R)).$$ 
Moreover, as $\vu \in L^{\infty}(]a,b[, L^{2}(B(x_{0},r)))$, for any $1\leq i\leq 3$, we have $(\partial_{i}\phi)\vu\in L^{\infty}(]0,+\infty[, L^{2}(\R))$, which results in 
$$\sum_{i=1}^{3}\partial_{i}\big((\partial_{i}\phi)\vu\big) \in L^{\infty}(]0, +\infty[, \dot{H}^{-1}(\R)).$$
With the two informations above, we get
$$\varphi\left( \frac{1}{\Delta} \left( -(\Delta \phi) \vu + 2\sum_{i=1}^{3}\partial_{i}\big((\partial_{i}\phi)\vu\big) \right) \right)\in L^{\infty}(]0, +\infty[, \dot{H}^{1}(\R)).$$ 
Hence, \eqref{claim} is verified by the Sobolev embedding $\dot{H}^{1}(\R) \subset L^{6}(\R)$. Once we have  $\vU-\vV \in L^{\infty}_t L_x^6$, by the assumption $5<q_0<6$ and by the localization property given in Lemma \ref{Lemme_locindi}, we have $\mathds{1}_{\Omega}(\vU-\vV) \in L_t^{q_0}L_x^{q_0}=M^{q_0, q_0}_{t,x} $ and this conclusion is enough for our purposes. However, let us note that, since $\sigma_0=\min\{\frac{p_{0}}{\lambda_{0}},q_{0}\} < q_0$, the fact  $\mathds{1}_{\Omega}(\vU-\vV)  \in  M_{t,x}^{\sigma_0,q_0}(\mathbb{R}\times\R)$ also follows from Lemma \ref{Lemme_locindi}. To finish, remark now that as we have the information $\vb \in L^{\infty}(]a,b[, L^{2}(B(x_{0},r)))$ and $\sigma_1=\min\{\frac{p_{1}}{\lambda_{1}},q_{1}\} < q_1<6$, the proof of the fact $\mathds{1}_{\Omega}(\vB-\vW) \in M_{t,x}^{\sigma_1,q_1}(\mathbb{R}\times\R)$ follows the same lines. \hfill$\blacksquare$\\

Once we have Proposition \ref{Propo_UmoinsV} for the differences $\vU-\vV$ and $\vB-\vW$, it remains to show that the quantities $\vV$ and $\vW$ defined in (\ref{Def_VW}) belong to the parabolic Morrey spaces $M_{t,x}^{\sigma_0,q_0}$ and $M_{t,x}^{\sigma_1,q_1}$. For this we will use the equations satisfied by these objects $\vV$ and $\vW$, but these dynamics involve the pressure $p$ for which we do not have any information  (recall that $p\in \mathcal{D}'$) and we need to get rid of this term, however, as we are working in a local setting we can not just apply the Leray projector and it will be more convenient to work with the \emph{vorticity} 
$$\vomega =  \vn \wedge \vu,$$ 
and with the \emph{current} 
$$\vrho = \vn \wedge \vb,$$
and with the equations satisfied by these two variables, which do not involve the pressure anymore: indeed if we apply the curl to the system (\ref{EcuacionMHD}) and since $\vn \wedge \vn p=0$ we will obtain the dynamics for $\vomega$ and $\vrho$ where there is no pressure.\\ 

The link between the variables $\vV, \vW$ defined in (\ref{Def_VW}) above and the functions $\vomega, \vrho$ is given by the following property: if we localize properly the vorticity $\vomega$ and the current $\vrho$, then due to the support properties of the localizing functions and by Lemma \ref{vort_vel} in the Appendix, we obtain (locally) the identities  
\begin{equation}\label{ProblemaVW}
\vV = - \varphi\left(\frac{1}{\Delta}(\phi\vn \wedge (\varphi \vomega))\right)=- \varphi\left(\frac{1}{\Delta}\phi\;\vUU\right)\quad \mbox{and} \quad \vW= -\varphi\left(\frac{1}{\Delta}(\phi\vn \wedge (\varphi\vrho))\right)=-\varphi\left(\frac{1}{\Delta}\phi\;\vBB\right),
\end{equation}
where $\vUU:=\vn\wedge (\varphi \vomega)$ and $\vBB:=\vn\wedge (\varphi \vrho)$. Thus in order to study $\vV$ and $\vW$ we shall first obtain some properties on the variables $\vUU$ and $\vBB$ since once we obtain information them it will be easy to deduce information for $\vV$ and $\vW$. Note that the dynamics for $\vUU$ and $\vBB$ can be deduced from the initial system (\ref{EcuacionMHD}) by first apply the curl, by localizing with the function $\varphi$ and by applying the curl again, we thus obtain the two following equations:
\begin{eqnarray}
\partial_{t}\vUU&=&\Delta \vUU\notag\\ 
&+&\vn\wedge\left[ \varphi(\vn \wedge \vf)+(\partial_{t}\varphi + \Delta \varphi ) \vomega - 2 \sum_{i=1}^{3}\partial_{i}((\partial_{i}\varphi)\vomega)  + \varphi\left(\vn \wedge\big(-(\vu\cdot \vn) \vu + (\vb\cdot \vn) \vb   \big)\right)\right]\label{DynamiqueUUBB}\\[3mm]  
\partial_{t}\vBB&=&\Delta \vBB\notag\\
& +&\vn\wedge\left[ \varphi(\vn \wedge \vg)+(\partial_{t}\varphi + \Delta \varphi ) \vrho - 2 \sum_{i=1}^{3}\partial_{i}((\partial_{i}\varphi)\vrho) 
+\varphi\left(\vn \wedge\big(-(\vu\cdot \vn) \vb +(\vb\cdot \vn) \vu \big)\right)\right]. \notag
\end{eqnarray}	
Remark now that by the definition of the localization function $\varphi$, we have $\vUU(0, \cdot)=0$ and $\vBB(0,\cdot)=0$ and thus the variables $\vUU$ and $\vBB$ satisfy the following parabolic equations: 
\begin{equation}\label{paraeqU}
\begin{cases}
\partial_{t}\vUU=\Delta\vUU+\vn\wedge \vec{\mathcal{R}} ,\\[2mm]
\vUU(0, \cdot)=0, 
\end{cases}
\qquad \mbox{and}\qquad 
\begin{cases}
\partial_{t}\vBB=\Delta\vBB+\vn\wedge \vec{\mathcal{V}},\\[2mm]
\vBB(0, \cdot)=0, 
\end{cases}
\end{equation}	
where, 
\begin{equation}\label{GrosseFormule1}
\begin{split}
\vec{\mathcal{R}}&= \sum_{j=1}^{11} \vec{\mathcal{R}}_j=\underbrace{\varphi(\vn \wedge \vf)}_{(1)}+\underbrace{(\partial_{t}\varphi + \Delta \varphi ) \vomega}_{(2)} - \underbrace{2 \sum_{i=1}^{3}\partial_i\big((\partial_{i}\varphi)\vomega\big)}_{(3)} \\ 
&+\underbrace{\sum_{i=1}^{3}  \partial_i (\vn \varphi \wedge (u_i \vu))}_{(4)} +\underbrace{ \vn \wedge \left( \sum_{i=1}^{3} (\partial_i \varphi) u_i \vu \right)}_{(5)}- \underbrace{\sum_{i=1}^{3}(\vn \partial_i \varphi) \wedge (u_i \vu)}_{(6)}- \underbrace{\vn \wedge\left(\sum_{i=1}^{3} \partial_i (\varphi u_i \vu)  \right)}_{(7)} \\
&- \underbrace{\sum_{i=1}^{3}  \partial_i (\vn \varphi \wedge (b_i \vb))}_{(8)} - \underbrace{\vn \wedge \left( \sum_{i=1}^{3} (\partial_i \varphi) b_i \vb \right)}_{(9)}+ \underbrace{\sum_{i=1}^{3}(\vn \partial_i \varphi) \wedge (b_i \vb)}_{(10)}+ \underbrace{\vn \wedge\left(\sum_{i=1}^{3} \partial_i (\varphi b_i \vb) \right)}_{(11)}
\end{split}
\end{equation}	
and 
\begin{equation*}
\begin{split}
\vec{\mathcal{V}}&= \sum_{j=1}^{11}\vec{\mathcal{V}}_j = \varphi(\vn \wedge \vg)+(\partial_{t}\varphi + \Delta \varphi ) \vrho - 2 \sum_{i=1}^{3}\partial_i\big((\partial_{i}\varphi)\vrho\big) \\ 
& +\sum_{i=1}^{3}  \partial_i (\vn \varphi \wedge (u_i \vb))+ \vn \wedge \left( \sum_{i=1}^{3} (\partial_i \varphi) u_i \vb \right)- \sum_{i=1}^{3}(\vn \partial_i \varphi) \wedge (u_i \vb) - \vn \wedge\left(\sum_{i=1}^{3} \partial_i (\varphi u_i \vb)  \right) \\
&-\sum_{i=1}^{3}  \partial_i (\vn \varphi \wedge (b_i \vu)) - \vn \wedge \left( \sum_{i=1}^{3} (\partial_i \varphi) b_i \vu \right)+ \sum_{i=1}^{3}(\vn \partial_i \varphi) \wedge (b_i \vu) + \vn \wedge\left(\sum_{i=1}^{3} \partial_i (\varphi b_i \vu)  \right).
\end{split}
\end{equation*}		
In the expressions of the quantities $\vec{\mathcal{R}}$ and $\vec{\mathcal{V}}$ given above we have systematically decomposed the terms $(\vu\cdot \vn)\vu$, $(\vb\cdot \vn)\vb$, $(\vu\cdot \vn)\vb$ and $(\vb\cdot \vn)\vu$ of (\ref{DynamiqueUUBB}) by using the identity
$$\varphi(\vn \wedge (\vA \cdot \vn) \vB ) =    - \sum_{i=1}^{3} \partial_i(\vn\varphi \wedge (A_i \vB)) - \vn \wedge (\sum_{i=1}^{n} (\partial_i \varphi) A_i \vB)   + \sum_{i=1}^{3} (\vn \partial_i \varphi) \wedge (A_i \vB) + \vn \wedge (\sum_{i=1}^{3} \partial_i(\varphi A_i \vB)),$$ 
which follows from vectorial identities, the support properties of the localization functions and the fact that $div(\vu)=div(\vb)=0$. See Lemma \ref{nonlinearident} in the Appendix for a detailed proof.\\
	
Now, using an integral representation we have that the solutions of equations \eqref{paraeqU} can be written in the following form
\begin{equation*}\label{eq06}
\vUU = \displaystyle\int_{0}^{t} e^{ (t-s)\Delta}  (\vn \wedge\vec{\mathcal{R}}) (s,\cdot)\,ds = \sum_{j=1}^{11}  \vn \wedge\displaystyle\int_{0}^{t} e^{ (t-s)\Delta}  \vec{\mathcal{R}}_{j} (s,\cdot)\,ds :=  \sum_{j=1}^{11} \vn \wedge\vec{T}_j,
\end{equation*}
and 
\begin{equation*}\label{eq07}
\vBB = \displaystyle\int_{0}^{t} e^{ (t-s)\Delta} (\vn \wedge \vec{\mathcal{V}}) (s,\cdot)\,ds=  \sum_{j=1}^{11} \vn \wedge\displaystyle\int_{0}^{t} e^{ (t-s)\Delta} \vec{\mathcal{V}}_{j} (s,\cdot)\,ds := \sum_{j=1}^{11}\vn \wedge\vec{X}_j,
\end{equation*}
where we defined $\vec{T}_j=\displaystyle\int_{0}^{t} e^{ (t-s)\Delta} \vec{\mathcal{R}}_{j} (s,\cdot)\,ds$ and $\vec{X}_j=\displaystyle\int_{0}^{t} e^{ (t-s)\Delta} \vec{\mathcal{V}}_{j} (s,\cdot)\,ds$.\\

With these expressions for the variables $\vUU$ and $\vBB$, we remark that in order to prove that $\vV\in M_{t,x}^{\sigma_0,q_0}$ and $\vW\in M_{t,x}^{\sigma_1,q_1}$, due to the identification (\ref{ProblemaVW}) we only have to verify that for each $\vec{T}_j$ and $\vec{X}_j$, with $j= 1,..., 11$, we actually have
\begin{equation}\label{goal1}
\varphi\Big(\frac{1}{\Delta}\big(\phi \vn \wedge\vec{T}_j\big)\Big) \in  M_{t,x}^{\sigma_0,q_0}(\mathbb{R}\times\R) \quad \text{and} \quad \varphi\Big(\frac{1}{\Delta}\big(\phi  \vn \wedge\vec{X}_j\big)\Big) \in  M_{t,x}^{\sigma_1,q_1}(\mathbb{R}\times\R).
\end{equation}
The rest of the paper is thus devoted to show \eqref{goal1} and for this we will treat separately each ones of the previous terms: indeed Proposition \ref{Propo_123} studies the cases $j=1,2$, Proposition \ref{Propo_3} treats the case $j=3$ while Proposition \ref{Propo_4to10} treats the cases $j=4,5,6,8,9,10$, finally Proposition \ref{Propo_7to11} studies the remaining cases, \emph{i.e.} $j=7,11$. 
\begin{Proposition}\label{Propo_123} Under the above notation, for $j=1,2$  we have 
$$\varphi\left(\frac{1}{\Delta} \big(\phi\vn \wedge \vec{T}_j\big)\right) \in M_{t,x}^{\sigma_0,q_0}\qquad \mbox{and}\qquad \varphi\left(\frac{1}{\Delta}\big(\phi \vn \wedge\vec{X}_j\big)\right) \in M_{t,x}^{\sigma_1,q_1}.$$
\end{Proposition}	
{\bf Proof of Proposition \ref{Propo_123}.} Let us start with $\vec{T}_1$. By Lemma \ref{Lemme_locindi}, since $5<q_0<6$ and since $\sigma_0=\min\{\frac{p_0}{\lambda_0}, q_0\}\leq q_0$ and using the identification $M_{t,x}^{p,p}=L^p_{t}L^p_x$, we can write
\begin{eqnarray*}
\left\|\varphi\left(\frac{1}{\Delta} \big(\phi \vn \wedge\vec{T}_1\big)\right)\right\|_{M_{t,x}^{\sigma_0,q_0}}&\leq &C\left\|\varphi\left(\frac{1}{\Delta} \big(\phi \vn \wedge\vec{T}_1\big)\right)\right\|_{L_{t}^{6}L_x^6}\leq C\|\varphi\|_{L^\infty_tL^\infty_x}\left\|\frac{1}{\Delta} \big(\phi \vn \wedge\vec{T}_1\big)\right\|_{L_{t}^{6}L_x^6}\notag\\
&\leq &C\left\|\frac{1}{\Delta} \big(\phi\vn \wedge \vec{T}_1\big)\right\|_{L_{t}^{6}\dot{H}_x^1}\leq C\|\phi \vn \wedge\vec{T}_1\|_{L_{t}^{6}\dot{H}_x^{-1}},\label{EstimateInter1}
\end{eqnarray*}
where we used the embedding $\dot{H}^1\subset L^6$ and the properties of the negative powers of the Laplacian. Now, by the definition of $\vec{T}_1$, using the embedding $ L^{\frac{6}{5}}\subset\dot{H}^{-1}$ and the H\"older inequality with $\tfrac{5}{6}=\tfrac{1}{3}+\tfrac{1}{2}$,  we write:
\begin{eqnarray*}
\|\phi \vn \wedge\vec{T}_1\|_{L_{t}^{6}\dot{H}_x^{-1}}&\leq &C\|\phi\vn \wedge \vec{T}_1\|_{L_{t}^{\infty}\dot{H}_x^{-1}}\leq C \|\phi\vn \wedge \vec{T}_1\|_{L_{t}^{\infty}L_x^{\frac{6}{5}}}= C\underset{t>0}{\sup}\left\|\phi\vn \wedge\int_0^t e^{(t-s)\Delta}\mathcal{R}_1ds\right\|_{L^{\frac{6}{5}}}\\
&\leq &C\|\phi\|_{L^\infty_tL^3_x}\;\underset{t>0}{\sup}\left\|\int_0^t e^{(t-s)\Delta}(-\Delta)^{\frac{1}{2}}\frac{(\vn \wedge \mathcal{R}_1)}{(-\Delta)^{\frac{1}{2}}}ds\right\|_{L^2}\\
&\leq &C\left\|\frac{\vn \wedge \mathcal{R}_1}{(-\Delta)^{\frac{1}{2}}}\right\|_{L_t^2L_x^{2}}= C\left\|\vn \wedge \mathcal{R}_1\right\|_{L_t^2\dot{H}_x^{-1}}.
\end{eqnarray*}
The last estimate follows from the general inequality $\underset{t>0}{\sup}\left\|\displaystyle{\int_0^t e^{(t-s)\Delta}}(-\Delta)^{\frac{1}{2}}F ds\right\|_{L^2}\leq C\|F\|_{L^2_tL^2_x}$, see Lemma \ref{Lema_MaximalRegularity} in the Appendix. Remark now that since $\vf\in L^{2}_{t}\dot{H}^{1}_{x}$ and due to the properties of the localizing function $\varphi$, we actually have $\vn \wedge \vec{\mathcal{R}}_{1}=\vn \wedge (\varphi(\vn \wedge \vf)) \in L^{2}_{t}\dot{H}^{-1}_{x} $ since
$$\|\vn \wedge \vec{\mathcal{R}}_{1}\|_{L^{2}_{t}\dot{H}^{-1}_{x}}\leq C\|\varphi(\vn \wedge \vf)\|_{L^{2}_{t}L^{2}_{x}}\leq C\|\vf\|_{L^{2}_{t}\dot{H}^{1}_{x}}<+\infty,$$ 
which finally gives $\varphi\left(\frac{1}{\Delta} \big(\phi \vn \wedge\vec{T}_1\big)\right) \in M_{t,x}^{\sigma_0,q_0}$. For $\vec{T}_2$, in a similar fashion, since we have by hypothesis $\vu\in L^{\infty}_{t}L^{2}_{x}\cap L^{2}_{t}\dot{H}^{1}_{x}(]a,b[\times B(x_0, r)) $ we get $\vn\wedge \vec{\mathcal{R}}_2= \vn\wedge\left((\partial_t\varphi+\Delta \varphi)(\vn\wedge \vu)\right)\in L_{t}^2 \dot{H}_{x}^{-1}$ from which we deduce that $\varphi\left(\frac{1}{\Delta} \big(\phi \vn \wedge\vec{T}_2\big)\right) \in M_{t,x}^{\sigma_0,q_0}$.\\

The estimates for $\varphi\left(\frac{1}{\Delta}\big(\phi \vn \wedge\vec{X}_j\big)\right)$ follow the same lines.  \hfill $\blacksquare$\\

\begin{Proposition}\label{Propo_3} For $j=3$,  we have 
$$\varphi\left(\frac{1}{\Delta}(\phi \vn \wedge \vec{T}_j) \right) \in M_{t,x}^{\sigma_0,q_0}\qquad \mbox{and} \qquad\varphi\left(\frac{1}{\Delta}(\phi \vn \wedge \vec{X}_j) \right) \in M_{t,x}^{\sigma_1,q_1}.$$
\end{Proposition}	
{\bf Proof of Proposition \ref{Propo_3}.}  
We will detail the first term since the second term that involves $\vec{X}_j$ follows the same computations.
Indeed, following the same ideas as previously we have
\begin{align*}
\left\|\varphi\left(\frac{1}{\Delta}(\phi \vn \wedge \vec{T}_3)\right) \right\|_{M_{t,x}^{\sigma_0,q_0}} & \leq C \left\|\varphi\left(\frac{1}{\Delta}(\phi \vn \wedge \vec{T}_3)\right) \right\|_{M_{t,x}^{q_0,q_0}} \leq C \left\|\varphi\left(\frac{1}{\Delta}(\phi \vn \wedge \vec{T}_3)\right) \right\|_{L^6_{t}L^6_x}.
\end{align*}
Let us define now $\vn \wedge \vec{T}_3 : = \Delta \vec{Y}_3$, where 
$$\vec{Y}_3 = -2 \sum_{i=1}^{3} \displaystyle\int_{0}^{t} e^{ (t-s)\Delta} \frac{1}{\Delta} \vn \wedge \partial_i \big( (\partial_{i}\varphi)\vomega\big) (s,\cdot)\,ds. $$
Using the classical identity 
$\phi (\Delta \vec{Y}_3)
= \Delta(\phi \vec{Y}_3)+(\Delta \phi)\vec{Y}_3 -2\displaystyle{\sum_{i=1}^{3}}\partial_{i}((\partial_{i}\phi)\vec{Y}_3)$, we obtain
\begin{eqnarray}
\varphi\Big(\frac{1}{\Delta}\big(\phi \vn \wedge \vec{T}_3\big) \Big)&=&\varphi\phi \vec{Y}_{3}+\varphi \frac{1}{\Delta}\left((\Delta \phi) \vec{Y}_{3}\right)-2 \sum_{i=1}^{3} \varphi \frac{\partial_{i}}{\Delta}\left(\left(\partial_{i} \phi\right) \vec{Y}_{3}\right).\label{term3deri}
\end{eqnarray}
It remains to treat each term on the right-hand side of equality \eqref{term3deri}. For the first term above, by using Sobolev embedding  $\dot{H}^{1}(\R) \subset L^{6}(\R)$ and a standard heat kernel estimate (see Lemma \ref{Lema_MaximalRegularity}), we get
\begin{eqnarray}
\|\varphi\phi \vec{Y}_{3} \|_{L^6_{t}L^6_x}	&\leq &C\| \vec{Y}_{3} \|_{L^\infty_{t} L^6_x} \leq C \| \vec{Y}_{3} \|_{L^\infty_{t} \dot H^1_x}\leq C\sum_{i=1}^{3} \left\| \int_{0}^{t} e^{ (t-s)\Delta} \left(\frac{1}{\Delta} \vn \wedge \partial_i \big( (\partial_{i}\varphi)\vomega\big)\right) (s,\cdot)\,ds\right\|_{L^\infty_{t} \dot H^1_x}\notag\\
&\leq &C \sum_{i=1}^{3}\left\| \frac{1}{\Delta} \vn \wedge \partial_i \big( (\partial_{i}\varphi)\vomega\big)\right\|_{L^2_{t} L^2_x}\leq C\sum_{i=1}^{3}\|(\partial_{i}\varphi)\vomega\|_{L^2_{t} L^2_x} \leq C \| \vu \|_{L^2_{t} \dot H^1_x} .\label{1term}
\end{eqnarray}
For the second term and the third term on the right-hand side of \eqref{term3deri}, we define the following two operators: 
$$
\mathcal{L}_1 : f \mapsto \varphi \frac{1}{\Delta}\left((\Delta \phi) f\right), \quad
\mathcal{L}_{2,i} : f \mapsto \varphi \frac{1}{\Delta} \partial_i \left((\partial_i \phi) f\right),
$$
which can be rewritten in the following form:
\begin{eqnarray*}
\mathcal{L}_1(f)(t,x)&=&\varphi \left(\frac{1}{\Delta}((\Delta\phi) f)\right)(t,x)=\varphi ( K\ast(\Delta \phi) f)(t,x)=C\varphi(t,x)\int_{\mathbb{R}^3}\frac{1}{|x-y|}\Delta\phi(t,y)f(t,y)dy\\
\mathcal{L}_{2,i}(f)(t,x)&=&\varphi \left(\frac{1}{\Delta}\partial_i((\partial_i\phi) f)\right)(t,x)=\varphi ( \partial_i K\ast(\partial_i \phi) f)(t,x)=C\varphi(t,x)\int_{\mathbb{R}^3}\frac{1}{|x-y|^2}\partial_i\phi(t,y)f(t,y)dy.
\end{eqnarray*}
By localization property of the test functions $\varphi$ and $\phi$, we find that $\mathcal{L}_1$ and $\mathcal{L}_{2,i}$ are bounded on $L^6_tL^6_x$, so we can write:
$$\left\|\varphi \frac{1}{\Delta}\left((\Delta \phi) \vec{Y}_{3}\right) \right\|_{L^6_{t}L^6_x}\leq C \| \mathds{1}_{ \textcolor{blue}{\Omega}}\vec{Y}_{3} \|_{L^6_{t}L^6_x}
 \leq C \| \vec{Y}_{3} \|_{L^\infty_{t} L^6_x},$$
and form the previous calculus displayed in (\ref{1term}) we finally obtain
$$ \left\|\varphi \frac{1}{\Delta}\left((\Delta \phi) \vec{Y}_{3}\right) \right\|_{L^6_{t}L^6_x}\leq C \| \vu \|_{L^2_{t} \dot H^1_x} <+\infty,$$
and we also have, for the last term of (\ref{term3deri})
\begin{align*}
\left\|\sum_{i=1}^{3} \varphi \frac{\partial_{i}}{\Delta}\left(\left(\partial_{i} \phi\right) \vec{Y}_{3}\right)\right\|_{L^6_{t}L^6_x}\leq \sum_{i=1}^{3} \left\|\mathcal{L}_{2,i} (\vec{Y}_{3}) \right\|_{L^6_{t}L^6_x}
& \leq C \| \mathds{1}_{\Omega_0} \vec{Y}_{3} \|_{L^6_{t}L^6_{x}} \leq C \| \vec{Y}_{3} \|_{L^\infty_{t} L^6_x}\leq C \| \vu \|_{L^2_{t} \dot H^1_x}. 
\end{align*}
Thus gathering all the $L^6_tL^6_x$ estimates for (\ref{term3deri}), we obtain $\varphi(\frac{1}{\Delta}(\phi (\vn \wedge \vec{T}_3)) ) \in M_{t,x}^{\sigma_0,q_0}$. \hfill$\blacksquare$\\

We continue our study of the terms $\vec{T_j}$ and $\vec{X_j}$ for $j=4,5,6,8,9,10$ and for this we will need to establish some estimates that involve the parabolic Riesz potencial $\mathcal{I}_\alpha$ defined in (\ref{Def_Riesz_potential}). 
\begin{Lemme}\label{obsvTvX} Under the notation above, for $j=4,5,6,8,9,10$ there exists a generic constant $C>0$ depending only on the size of the set $\Omega=]a,b[\times B(x_{0},r)$, such that the variables $\vec{T_j}$ and $\vec{X_j}$ verify the following pointwise estimates: 
\begin{enumerate}
\item[1)] For $j=4,5$,  \,\,   
$\ds{\vert \vec{T_j}(t,x) \vert \leq C \, \mathcal{I}_{1}(\mathds{1}_{\Omega}  \vert \vu(t,x)\vert^2) }$ and $\ds{ \vert \vec{X_j}(t,x) \vert \leq C\, \mathcal{I}_{1}(\mathds{1}_{\Omega}  \vert \vu(t,x) \otimes \vb(t,x)\vert) }$.  \\
\item[2)] For $j=6$, \, \, $\ds{\vert \vec{T_j}(t,x) \vert \leq C \, \mathcal{I}_{2}(\mathds{1}_{\Omega}  \vert \vu(t,x)\vert^2) }$ and   $ \ds{ \vert \vec{X_j}(t,x) \vert \leq C \, \mathcal{I}_{2}(\mathds{1}_{\Omega}  \vert \vu(t,x) \otimes \vb(t,x)\vert) }$.  \\
\item[3)] For $j=8,9$, \,\,  $\ds{\vert \vec{T_j}(t,x) \vert \leq C  \, \mathcal{I}_{1}(\mathds{1}_{\Omega}  \vert \vb(t,x)\vert^2) }$ and  $ \ds{ \vert \vec{X_j}(t,x) \vert \leq C \, \mathcal{I}_{1}(\mathds{1}_{\Omega}  \vert \vu(t,x) \otimes \vb(t,x)\vert) }$.  \\
\item[4)] For $j=10$, \,\,  $\ds{\vert \vec{T_j}(t,x) \vert \leq C \, \mathcal{I}_{2}(\mathds{1}_{\Omega}  \vert \vb(t,x)\vert^2) }$ and  $ \ds{ \vert \vec{X_j}(t,x) \vert \leq C \, \mathcal{I}_{2}(\mathds{1}_{\Omega}  \vert \vu(t,x) \otimes \vb(t,x)\vert) }$.  \\ 
\end{enumerate}		
\end{Lemme}	
{\bf Proof of Lemma \ref{obsvTvX}.}  We detail here only the estimates for the values $j=4$ and $j=6$ since the proofs of all the other terms follow essentially the same computations due to their common structure. 
\begin{enumerate}
\item[$\bullet$] For $j=4$, recalling that we have the following expression for the heat semi-group $e^{ (t-s)\Delta}f=h_{t-s}\ast f$ where $h_t$ is the heat kernel, we can write 
\begin{eqnarray*}
\vec{T_4}(t,x)&=&\int_{0}^{t} e^{ (t-s)\Delta}  \left( \sum_{i=1}^{3}  \partial_i (\vn \varphi \wedge (u_i \vu) ) \right) (s,x)ds=\sum_{i=1}^{3}\int_{0}^{t}  \int_{\mathbb{R}^3}  \partial_i h_{t-s}(x-y) \vn \varphi \wedge (u_i \vu) (s,y)dyds\\
|\vec{T_4}(t,x)|&\leq &\sum_{i=1}^{3}\int_{0}^{t}\int_{\mathbb{R}^3}  \left| \partial_i h_{t-s}(x-y) \right| |\vn \varphi \wedge (u_i \vu)| (s,y)dyds.
\end{eqnarray*}
By the decay properties of the heat kernel (see Lemma \ref{Lema_Decay} in the Appendix) and by the support properties of the function $\varphi$, we observe that we have
$$\vert \vec{T_4}(t,x) \vert \leq C  \sum_{i=1}^{3}  \displaystyle\int_{\mathbb{R}} \displaystyle\int_{\mathbb{R}^3} \frac{1}{(|t-s|^{\frac{1}{2}}+|x-y|)^4}|\vn \varphi \wedge (u_i \vu)(s,y)|\,dy \,ds,$$
now, with the definition of the parabolic Riesz potential $\mathcal{I}_1$ given in (\ref{Def_Riesz_potential}) and with the boundedness properties of the function $\varphi$ we have:
$$|\vec{T_4}(t,x)| \leq  C \sum_{i=1}^{3}  \, \mathcal{I}_{1}(|\vn \varphi \wedge (u_i \vu)|) \leq    C \, \mathcal{I}_{1}(\mathds{1}_{\Omega}  \vert \vu(t,x)\vert^2). $$
Remember that $\vec{X_4}(t,x)$ has the same expression as $\vec{T_4}(t,x)$ by replacing $\vu$ by $\vb$. So we may use the same technique to show that 
$$ \ds{ \vert \vec{X_4}(t,x) \vert \leq C \, \mathcal{I}_{1}(\mathds{1}_{\Omega}  \vert \vu(t,x) \otimes \vb(t,x)\vert) }.$$
\item[$\bullet$] For $j=6$, recall that we have 
\begin{equation*}
\begin{split}
\vec{T_6}(t,x)= \int_{0}^{t} e^{ (t-s)\Delta}  \left(\sum_{i=1}^{3}(\vn \partial_i \varphi) \wedge (u_i \vu) \right)  (s,x) ds,
\end{split}
\end{equation*}
and by the same arguments above we can write
\begin{eqnarray*}
\vert \vec{T_6}(t,x) \vert 
& \leq  & C \sum_{i=1}^{3}  \displaystyle\int_{\mathbb{R}}
\displaystyle\int_{\mathbb{R}^3}\frac{1}{(|t-s|^{\frac{1}{2}}+|x-y|)^3}|\vn \partial_i \varphi \wedge (u_i \vu) (s,y)|  \,dy \,ds  \leq C \, \mathcal{I}_{2}(\mathds{1}_{\Omega}  \vert \vu(t,x)\vert^2). 
\end{eqnarray*}
\end{enumerate}
The same computations for $\vec{X_6}(t,x)$ lead us to obtain $ \ds{ \vert \vec{X_6}(t,x) \vert \leq C\, \mathcal{I}_{2}(\mathds{1}_{\Omega}  \vert \vu(t,x) \otimes \vb(t,x)\vert) }$.  \hfill$\blacksquare$\\

Once we have these pointwise estimates, we can continue our study where we will use the hypothesis on $\vu$ and $\vb$ given in (\ref{LocalHypo1}).  
\begin{Proposition} \label{Propo_4to10} Under the notation above and for $j=4,5,6,8,9,10$ we have 
$$\varphi\left(\frac{1}{\Delta}(\phi \vn \wedge \vec{T}_j) \right) \in M_{t,x}^{\sigma_0,q_0}\qquad  \mbox{and}\qquad \varphi\left(\frac{1}{\Delta}(\phi \vn \wedge \vec{X}_j) \right) \in M_{t,x}^{\sigma_1,q_1}.$$
\end{Proposition}
{\bf Proof of Proposition \ref{Propo_4to10}}
As for the previous lemma, we will only detail here some cases since the proof of the remaining cases follows essentially the same computations.  
\begin{enumerate}
\item[$\bullet$] For the term $\varphi\Big(\frac{1}{\Delta}\big(\phi \vn \wedge \vec{T}_4\big) \Big)$, we have 
\begin{eqnarray}
\left\|\varphi\Big(\frac{1}{\Delta}\big(\phi \vn \wedge \vec{T}_4\big) \Big) \right\|_{M_{t,x}^{\sigma_0,q_0}}& = &\left\|\varphi\Big(\frac{1}{\Delta}\big(\vn \wedge(\phi \vec{T}_{4})-\vec{\nabla} \phi \wedge \vec{T}_{4}\big) \Big) \right\|_{M_{t,x}^{\sigma_0,q_0}}	\notag\\
& \leq &\left\|\varphi\Big(\frac{1}{\Delta}\vn \wedge(\phi \vec{T}_{4}) \Big) \right\|_{M_{t,x}^{\sigma_0,q_0}} + \left\|\varphi\Big(\frac{1}{\Delta}\vec{\nabla} \phi \wedge \vec{T}_{4} \Big) \right\|_{M_{t,x}^{\sigma_0,q_0}}.\label{EstimationT5}
\end{eqnarray} 
Let us remark now that the inner structure of the terms $\varphi\Big(\frac{1}{\Delta}\vn \wedge(\phi \vec{T}_{4}) \Big) $ and $\varphi\Big(\frac{1}{\Delta}\vec{\nabla} \phi \wedge \vec{T}_{4} \Big)$ is of the following form:
\begin{eqnarray*}
\mathcal{T}_{a,i}(f)(t,x)=\varphi \left(\frac{1}{\Delta}\partial_i(\phi f)\right)(t,x)&=&\varphi (\partial_i K\ast \phi f)(t,x)=C\varphi(t,x)\int_{\mathbb{R}^3}\frac{1}{|x-y|^2}\phi(t,y)f(t,y)dy\quad \mbox{and}\\
\mathcal{T}_{b,i}(f)(t,x)=\varphi \left(\frac{1}{\Delta}(\partial_i\phi) f\right)(t,x)&=&C\varphi(t,x)\int_{\mathbb{R}^3}\frac{1}{|x-y|}\partial_i\phi(t,y)f(t,y)dy,
\end{eqnarray*}
where $K$ is the convolution kernel associated to the operator $\frac{1}{\Delta}$ (namely $\frac{1}{|x|}$) and $f(t,x)$ is a component of the vector $\vec{T}_{4}$. At this point we observe that, due to the support properties of the functions $\varphi$ and $\phi$, the kernels associated to the operators $\mathcal{T}_{a,i}$ and $\mathcal{T}_{b,i}$ are bounded in $L^1(\mathbb{R}^3)$.
Since the norm of $M^{\sigma_{0}, q_{0}}$ is translation invariant, we deduce from these facts the estimates
$\|\mathcal{T}_{a,i}(f)\|_{M_{t,x}^{\sigma_{0}, q_{0}}}\leq C_a\|f\|_{M_{t,x}^{\sigma_{0}, q_{0}}}$ and $\|\mathcal{T}_{b,i}(f)\|_{M_{t,x}^{\sigma_{0}, q_{0}}}\leq C_b\|f\|_{M_{t,x}^{\sigma_{0}, q_{0}}}$. Applying these observations to the right-hand side of (\ref{EstimationT5}), and keeping in mind the support of the function $\varphi$ given in (\ref{DefSoporteFuncTest2}) we have
\begin{equation}\label{PetiteFormuleUtile1}
\left\|\varphi\Big(\frac{1}{\Delta}\big(\phi \vn \wedge \vec{T}_4\big) \Big) \right\|_{M_{t,x}^{\sigma_0,q_0}}\leq C_a \|\mathds{1}_{\Omega_1}\vec{T}_4\|_{M_{t,x}^{\sigma_0,q_0}}+C_b \|\mathds{1}_{\Omega_1}\vec{T}_4\|_{M_{t,x}^{\sigma_0,q_0}}.
\end{equation}
Now, using the first point of Lemma \ref{obsvTvX}, the second point of Corollary \ref{Coro_I1} and Lemma \ref{Lemme_Product} we obtain:
\begin{eqnarray*}
\left\|\varphi\Big(\frac{1}{\Delta}\big(\phi \vn \wedge \vec{T}_4\big) \Big) \right\|_{M_{t,x}^{\sigma_0,q_0}}&\leq &C\|\mathds{1}_{\Omega_1}  \mathcal{I}_{1}(\mathds{1}_{\Omega_0}  |\vu|^2)\|_{M_{t,x}^{\sigma_{0}, q_{0}}}\leq C\|\mathds{1}_{\Omega_0}  |\vu|^2 \|_{M_{t,x}^{\frac{p_{0}}{2}, \frac{q_{0}}{2}}}\\
&\leq & C \|\mathds{1}_{\Omega}  |\vu|\|_{M_{t,x}^{p_{0}, q_0}}^2,
\end{eqnarray*}
thus, by the assumption $\mathds{1}_{\Omega}\vu\in M_{t,x}^{p_{0},q_{0}}(\mathbb{R}\times\R)$, we can conclude that $\varphi\Big(\frac{1}{\Delta}\big(\phi \vn \wedge \vec{T}_4\big) \Big) \in M_{t,x}^{\sigma_0,q_0} (\mathbb{R}\times\R)$. 

\item[$\bullet$] For the term $\varphi\Big(\frac{1}{\Delta}\big(\phi \vn \wedge \vec{X}_4\big) \Big)$, we can perform the same computations as before to obtain
\begin{eqnarray*}
\left\|\varphi\Big(\frac{1}{\Delta}\big(\phi \vn \wedge \vec{X}_4\big) \Big) \right\|_{M_{t,x}^{\sigma_1,q_1}}&\leq &\left\|\varphi\Big(\frac{1}{\Delta}\vn \wedge(\phi \vec{X}_{4})\Big)\right\|_{M_{t,x}^{\sigma_1,q_1}}+\left\|\varphi\Big(\frac{1}{\Delta}\vec{\nabla} \phi \wedge \vec{X}_{4} \Big) \right\|_{M_{t,x}^{\sigma_1,q_1}}	\\
&\leq & C_a\left\|\mathds{1}_{\Omega_1} \vec{X}_{4}\right\|_{M_{t,x}^{\sigma_1,q_1}}+C_b\left\|\mathds{1}_{\Omega_1} \vec{X}_{4}\right\|_{M_{t,x}^{\sigma_1,q_1}}\\
&\leq & C\|\mathds{1}_{\Omega_1} \mathcal{I}_{1}(\mathds{1}_{\Omega}  \vert \vu(t,x) \otimes \vb(t,x)\vert) \|_{M_{t,x}^{\sigma_1,q_1}}\\
& \leq & C\|\mathds{1}_{\Omega_0}  \vert \vu \otimes \vb\vert \|_{M_{t,x}^{\frac{p_1}{2}, \frac{q_1}{2}}} \leq \|\mathds{1}_{\Omega}  \vert \vu\vert \|_{M_{t,x}^{p_1, q_1}} \|\mathds{1}_{\Omega}  \vert \vb\vert \|_{M_{t,x}^{p_1, q_1}} 
\end{eqnarray*} 
 By the hypotheses $\mathds{1}_{\Omega}\vu\in M_{t,x}^{p_{0},q_{0}}(\mathbb{R}\times\R)$ and $\mathds{1}_{\Omega}\vb\in M_{t,x}^{p_{1},q_{1}}(\mathbb{R}\times\R)$ with $p_1 \leq p_0$, $q_1 \leq q_0$, we finally obtain that $\varphi\Big(\frac{1}{\Delta}\big(\phi \vn \wedge \vec{X}_4\big) \Big)\in M_{t,x}^{p_{1},q_{1}}(\mathbb{R}\times\R) $ .\\

\end{enumerate}
For the cases $j=6$ and $j=10$ we can apply the same arguments, the only modification is given by the use of Corollary \ref{Coro_I2} in order to study the parabolic Riesz potential $\mathcal{I}_2$. \hfill$\blacksquare$\\
\begin{Proposition}\label{Propo_7to11}
For the remaining terms of (\ref{GrosseFormule1}), \emph{i.e.} for $j=7, 11$, we have 
$$\varphi\left(\frac{1}{\Delta}\big(\phi \vn \wedge \vec{T}_j\big) \right) \in M_{t,x}^{\sigma_0,q_0}\qquad \mbox{and}\qquad \varphi\left(\frac{1}{\Delta}\big(\phi \vn \wedge \vec{X}_j\big) \right) \in M_{t,x}^{\sigma_1,q_1}.$$
\end{Proposition} 
{\bf Proof of Proposition \ref{Propo_7to11}}	
We will detail the case $j=7$ as the case when $j=11$ follows the same computations. 
\begin{enumerate}
\item[$\bullet$] Recall that 
$$ \vn \wedge \vec{T}_{7} = -\int_{0}^{t}  \vn \wedge \vn \wedge e^{ (t-s)\Delta}   \Big(\sum_{i=1}^{3} \partial_i (\varphi u_i \vu)  \Big) \,ds$$
Let $\vn \wedge \vec{T}_{7} :=\Delta \vec{Y}_{7}$, more precisely we have
$$\vec{Y}_{7}=-\sum_{i=1}^{3} \int_{0}^{t} \frac{1}{\Delta} \vn \wedge \vn \wedge e^{ (t-s)\Delta}    \partial_i (\varphi u_i\vu) \, ds.$$
Using the classical identity $\phi (\Delta \vec{Y}_7)= \Delta(\phi \vec{Y}_7)+(\Delta \phi)\vec{Y}_7 -2\displaystyle{\sum_{i=1}^{3}\partial_{i}((\partial_{i}\phi)\vec{Y}_7)}$,
we obtain 
$$\varphi\Big(\frac{1}{\Delta}\big(\phi \vn \wedge \vec{T}_7\big) \Big)=\varphi\phi \vec{Y}_{7}+\varphi \frac{1}{\Delta}\left((\Delta \phi) \vec{Y}_{7}\right)-2 \sum_{i=1}^{3} \varphi \frac{\partial_{i}}{\Delta}\left(\left(\partial_{i} \phi\right) \vec{Y}_{7}\right),$$
which gives 
\begin{eqnarray}
\left\|\varphi\Big(\frac{1}{\Delta}\big(\phi \vn \wedge \vec{T}_7\big) \Big)\right\|_{M_{t,x}^{\sigma_0,q_0}}&\leq &\left\|\varphi\phi \vec{Y}_{7}\right\|_{M_{t,x}^{\sigma_0,q_0}}+\left\|\varphi \frac{1}{\Delta}\left((\Delta \phi) \vec{Y}_{7}\right)\right\|_{M_{t,x}^{\sigma_0,q_0}}\label{PetiteFormuleUtile2}\\
&&+2 \sum_{i=1}^{3} \left\|\varphi \frac{\partial_{i}}{\Delta}\left(\left(\partial_{i} \phi\right) \vec{Y}_{7}\right)\right\|_{M_{t,x}^{\sigma_0,q_0}},\notag
\end{eqnarray}
and we will study each one of the previous terms separately. For the first term above we write
\begin{eqnarray*}
\|\varphi\phi \vec{Y}_{7} \|_{M_{t,x}^{\sigma_0,q_0}}&\leq &\sum_{i=1}^{3} \left\|\varphi\phi\int_{0}^{t} \frac{1}{\Delta} \vn \wedge \vn \wedge e^{ (t-s)\Delta}    \partial_i (\varphi u_i\vu) \, ds\right \|_{M_{t,x}^{\sigma_0,q_0}}\\
&\leq &\sum_{i=1}^{3} \left\|\varphi\phi\int_{0}^{t} \partial_i   e^{ (t-s)\Delta}    \left(\frac{1}{\Delta}\vn \wedge \vn \wedge(\varphi u_i\vu)\right) \, ds\right \|_{M_{t,x}^{\sigma_0,q_0}}.
\end{eqnarray*}
Let us study the quantity $\varphi\phi\displaystyle{\int_{0}^{t}} \partial_i   e^{ (t-s)\Delta}  \left(\frac{1}{\Delta}\vn \wedge \vn \wedge(\varphi u_i\vu)\right) \, ds$ and we can write
$$\varphi\phi\displaystyle{\int_{0}^{t}} \partial_i   e^{ (t-s)\Delta}  \left(\frac{1}{\Delta}\vn \wedge \vn \wedge(\varphi u_i\vu)\right) \, ds=\varphi\phi\displaystyle{\int_{0}^{t}} \int_{\mathbb{R}^3}\partial_i   h_{t-s}(x-y) \left(\frac{1}{\Delta}\vn \wedge \vn \wedge(\varphi u_i\vu)\right)(s,y) \,dy ds,$$
and due to the support properties of the functions $\varphi$ and $\phi$, and to the decay properties of the heat kernel (see Lemma \ref{Lema_Decay} in the Appendix), we obtain
$$\left|\varphi\phi\displaystyle{\int_{0}^{t}} \partial_i   e^{ (t-s)\Delta}  \left(\frac{1}{\Delta}\vn \wedge \vn \wedge(\varphi u_i\vu)\right) \, ds\right|\leq|\varphi\phi|\int_{\mathbb{R}} \int_{\mathbb{R}^3}|\partial_i   h_{t-s}(x-y)| \left|\frac{1}{\Delta}\vn \wedge \vn \wedge(\varphi u_i\vu)(s,y)\right| \,dy ds$$
\begin{eqnarray*}
&\leq &C|\varphi\phi|\int_{\mathbb{R}} \int_{\mathbb{R}^3} \frac{1}{(|t-s|^{\frac{1}{2}}+|x-y|)^4}   \left|\frac{1}{\Delta}\vn \wedge \vn \wedge(\varphi u_i\vu)(s,y)\right| \,dy ds\\
&\leq &C|\varphi\phi| \mathcal{I}_1\left(\left|\frac{1}{\Delta}\vn \wedge \vn \wedge(\varphi u_i\vu)\right|\right), 
\end{eqnarray*}
and with this estimate we have
$$\sum_{i=1}^{3} \left\|\varphi\phi\int_{0}^{t} \partial_i   e^{ (t-s)\Delta}    \left(\frac{1}{\Delta}\vn \wedge \vn \wedge(\varphi u_i\vu)\right) \, ds\right \|_{M_{t,x}^{\sigma_0,q_0}}\leq C\sum_{i=1}^{3} \left\|\mathds{1}_{\Omega_0} \mathcal{I}_1\left(\left|\frac{1}{\Delta}\vn \wedge \vn \wedge(\varphi u_i\vu)\right|\right)\right\|_{M_{t,x}^{\sigma_0,q_0}}.$$
By the localization properties of $\phi$, the second point of Corollary \ref{Coro_I1} and the boundedness of Riesz transforms in Morrey spaces and Lemma \ref{Lemme_Product}, we have 
\begin{equation}\label{PetiteFormuleUtile3}
\|\varphi \phi \vec{Y}_{7} \|_{M_{t,x}^{\sigma_0,q_0}} \leq C\sum_{i=1}^{3} \left\| \frac{1}{\Delta}\vn \wedge \vn \wedge(\varphi u_i\vu)\right\|_{M_{t,x}^{\frac{p_{0}}{2}, \frac{q_{0}}{2}}}\leq C\|\mathds{1}_{\Omega_0}  \vert \vu(t,x)\vert \|^2_{M_{t,x}^{p_0, q_0}}<+\infty,
\end{equation}
since by hypothesis we have $\mathds{1}_{\Omega_0}\vu\in M_{t,x}^{p_0, q_0}$ and we thus obtain the wished estimate for the first term of the right-hand side of (\ref{PetiteFormuleUtile2}).\\

For the second and the third term of the right-hand side of (\ref{PetiteFormuleUtile2}), using the same strategy as in the proof of Proposition \ref{Propo_4to10} (see \eqref{EstimationT5} and \eqref{PetiteFormuleUtile1}) and with the previous estimate (\ref{PetiteFormuleUtile3}) we finally obtain        
\begin{eqnarray} 
\|\varphi \frac{1}{\Delta}\left((\Delta \phi) \vec{Y}_{7}\right) \|_{M_{t,x}^{\sigma_0,q_0}} \leq C \| \mathds{1}_{\Omega_1} \vec{Y}_{7} \|_{M_{t,x}^{\sigma_0,q_0}} \leq  C\|\mathds{1}_{\Omega_0}  \vert \vu(t,x)\vert \|^2_{M_{t,x}^{p_0, q_0}}<+\infty,
\label{EstimationY7_2}
\end{eqnarray} 
and
\begin{eqnarray} 
\|\sum_{i=1}^{3} \varphi \frac{\partial_{i}}{\Delta}\left(\left(\partial_{i} \phi\right) \vec{Y}_{7}\right) \|_{M_{t,x}^{\sigma_0,q_0}}
\leq C \| \mathds{1}_{\Omega_1} \vec{Y}_{7} \|_{M_{t,x}^{\sigma_0,q_0}}
\leq  C\|\mathds{1}_{\Omega_0}  \vert \vu(t,x)\vert \|^2_{M_{t,x}^{p_0, q_0}}<+\infty.
\label{EstimationY7_3}
\end{eqnarray} 
Gathering the relations \eqref{PetiteFormuleUtile3}-\eqref{EstimationY7_3}, we can conclude that each term of (\ref{PetiteFormuleUtile2}) is bounded and we have
$$\varphi\left(\frac{1}{\Delta}(\phi \vn \wedge \vec{T}_7) \right) \in M_{t,x}^{\sigma_0,q_0} (\mathbb{R}\times\R).$$ 
\item[$\bullet$] Recall that
$$ 
\vn \wedge \vec{X}_{7} = -\int_{0}^{t}  \vn \wedge \vn \wedge e^{ (t-s)\Delta}   \Big(\sum_{i=1}^{3} \partial_i (\varphi u_i \vb)  \Big) \,ds$$
Let us define $\vn \wedge \vec{X}_{7} :=\Delta \vec{Z}_{7} $.
As $\mathds{1}_{\Omega}\vu\in M_{t,x}^{p_{0},q_{0}}(\mathbb{R}\times\R)$ and $\mathds{1}_{\Omega}\vb\in M_{t,x}^{p_{1},q_{1}}(\mathbb{R}\times\R)$ with $p_0 \leq q_0$, $p_1 \leq q_1$, $p_1 \leq p_0$ and $q_1\leq q_0$ the same calculus can be used to complete our proof.  \hfill$\blacksquare$
\end{enumerate}

\section{Appendix}
\begin{Lemme}\label{Lema_Decay} Let $h_t$ be the heat kernel. If $\alpha\in \mathbb{N}^3$ is a multi-index then we have
$$\left|D^\alpha h_{t}(x)\right|\leq C 
\begin{cases}
|x|^{-(3+|\alpha|)} \mbox{ if } |x|^2> t,\\[3mm]
t^{-\frac{(3+|\alpha|)}{2}} \mbox{ if } |x|^2\leq t.
\end{cases}
$$
\end{Lemme}
See \cite{Saka} for a proof of these facts in a general framework.
\begin{Lemme}\label{Lema_MaximalRegularity}
If $f\in L^{2}([0,+\infty[, L^{2}(\mathbb{R}^{3}))$ and if we define $F(t,x)=\displaystyle{\int_{0}^{t}h_{t-s}\ast f(s,x)ds}$ then we have
$$\|F(t,\cdot)\|_{\dot{H}^{1}_{x}}\leq C\|f\|_{L^{2}_{t}L^{2}_{x}}.$$
\end{Lemme}
{\bf Proof of Lemma \ref{Lema_MaximalRegularity}.} We simply write
\begin{eqnarray*}
\|(-\Delta)^{\frac 12}F(t,\cdot)\|_{L^{2}}&=&\underset{\|\phi\|_{L^{2}}\leq 1}{\sup}\left|\int_{\R}(-\Delta)^{\frac 12} F(t,x) \phi(x)dx\right|=\underset{\|\phi\|_{L^{2}}\leq 1}{\sup}\left|\int_{\R} \int_{0}^{t} (-\Delta)^{\frac 12}\left(h_{t-s}\ast f(s,x)\right)ds \phi(x)dx\right|\\
&=&\underset{\|\phi\|_{L^{2}}\leq 1}{\sup}\left|\int_{0}^{t}\int_{\R}(-\Delta)^{\frac 12}\left(h_{t-s} \ast \phi\right)  f(s,x)dxds\right|\\
&\leq &\underset{\|\phi\|_{L^{2}}\leq 1}{\sup}\int_{0}^{t}\|f(s,\cdot)\|_{L^{2}}\|(-\Delta)^{\frac 12} (h_{t-s}\ast \phi)\|_{L^{2}}ds\leq  \underset{\|\phi\|_{L^{2}}\leq 1}{\sup}\|f\|_{L^{2}_{t}L^{2}_{x}}\| h_{t}\ast \phi\|_{L^{2}_{t}\dot{H}^{1}_{x}}.
\end{eqnarray*}
Now remark that we have for the last term above
$$\|h_{t}\ast \phi\|_{L^{2}_{t}\dot{H}^{1}_{x}}^{2}\simeq\int_{0}^{+\infty}\int_{\mathbb{R}^{3}}|\xi|^{2}e^{-2t|\xi|^{2}}|\widehat{\phi}(\xi)|^{2}d\xi dt=\int_{\mathbb{R}^{3}}\int_{0}^{+\infty}|\xi|^{2}e^{-2t|\xi|^{2}}|\widehat{\phi}(\xi)|^{2}dtd\xi,$$
thus, by the change of variable $\tau=2t|\xi|^{2}$ we can write $\displaystyle{\|h_{t}\ast \phi\|_{L^{2}_{t}\dot{H}^{1}_{x}}^{2}\simeq \int_{\mathbb{R}^{3}}\int_{0}^{+\infty} e^{-\tau}|\widehat{\phi}(\xi)|^{2}d\tau d\xi=\|\widehat{\phi}\|_{L^{2}}^2}$  
which gives $\|h_{t}\ast \phi\|_{L^{2}_{t}\dot{H}^{1}_{x}}\leq C\|\phi\|_{L^{2}}$, and we finally obtain $\|F(t,\cdot)\|_{\dot{H}^{1}_{x}}\leq C\|f\|_{L^{2}_{t}L^{2}_{x}}$. \hfill$\blacksquare$\\
\begin{Lemme}\label{vort_vel}
If $\phi$ is the test function defined in (\ref{DefSoporteFuncTest1}), if $\varphi$ the test function defined in (\ref{DefSoporteFuncTest2}) and $\vu$ is a regular vector field, then we have
$$\varphi\left(\frac{1}{\Delta}(\phi(\Delta \vu))\right)=-\varphi\left(\frac{1}{\Delta}\left(\phi \vn \wedge [\varphi \vn\wedge \vu]\right)\right).$$ 
\end{Lemme}
{\bf Proof of Lemma \ref{vort_vel}.} Indeed, 
\begin{eqnarray*}
\vn \wedge [\varphi \vn\wedge \vu]&=&\varphi \vn \wedge (\vn \wedge \vu)+\vn \varphi \wedge (\vn\wedge \vu)=\varphi \left(\vn(div(\vu))-\Delta \vu\right)+\vn \varphi \wedge (\vn\wedge \vu)\\
&=&-\varphi \Delta \vu+\vn \varphi \wedge (\vn\wedge \vu).
\end{eqnarray*}
Moreover, by the support properties of $\phi$ and $\varphi$ we have $\vn\varphi\equiv 0$ and $\phi\varphi=1$ on the support of $\phi$. So the second term in the identity above will disappear when we multiply the identity by $\phi$ and then we have 
$$\phi(\vn \wedge [\varphi \vn\wedge \vu])=\phi(-\varphi \Delta \vu+\vn \varphi \wedge (\vn\wedge \vu))=-\Delta \vu.$$
\hfill$\blacksquare$\\
For the nonlinear terms in the equations, we use the following lemma:
\begin{Lemme}\label{nonlinearident}
Let $\vA=(A_1, A_2, A_3)$ and $\vB=(B_1, B_2, B_3)$ be two functions such that $div(\vA)=0$ and $div(\vB)=0$. Then, 
\begin{equation}\label{eq01}
\begin{split}
\varphi(\vn \wedge (\vA \cdot \vn) \vB ) = &   - \sum_{i=1}^{3} \partial_i(\vn\varphi \wedge (A_i \vB)) - \vn \wedge (\sum_{i=1}^{n} (\partial_i \varphi) A_i \vB) \\ 
&  + \sum_{i=1}^{3} (\vn \partial_i \varphi) \wedge (A_i \vB) + \vn \wedge (\sum_{i=1}^{3} \partial_i(\varphi A_i \vB)).  
\end{split}
\end{equation}
\end{Lemme}
{\bf Proof.}  We write  $\varphi(\vn \wedge (\vA \cdot \vn) \vB )= \vn \wedge (\varphi(\vA \cdot \vn) \vB )- \vn \varphi \wedge ((\vA \cdot \vn) \vB ) $ where we study each term in the right-hand side. As $div(\vA)=0$ we can write 
\begin{eqnarray*}
\vn \wedge (\varphi(\vA \cdot \vn) \vB ) &=& \vn \wedge (\varphi \sum_{i=1}^{3} A_i \partial_i \vB ) = \vn \wedge ( \sum_{i=1}^{3} \varphi \partial_i( A_i  \vB ))	 \\
&=& \vn \wedge (\sum_{i=1}^{3} (\partial_i(\varphi A_i \vB) - (\partial_i \varphi)(A_i \vB))) \\
&=&  - \vn \wedge ((\partial_i \varphi) A_i \vB ) + \vn \wedge (\sum_{i=1}^{3} \partial_i(\varphi A_i \vB) ) , 
\end{eqnarray*} 
where we obtain the second and fourth terms in (\ref{eq01}). Then, always as we have $div(\vA)=0$ we write 
\begin{eqnarray*}
-\vn \varphi \wedge ((\vA \cdot \vn) \vB )&=& -\vn \varphi \wedge (\sum_{i=1}^{3} \partial_i (A_i \vB)) = - \sum_{i=1}^{3} \vn \varphi \wedge \partial_i(A_i \vB)  \\
&=& -\sum_{i=1}^{3} ( \partial_i(\vn \varphi \wedge (A_i \vB)) +\vn \partial_i \varphi \wedge (A_i \vB))\\
&=&- \sum_{i=1}^{3}  \partial_i(\vn \varphi \wedge (A_i \vB)) + \sum_{i=1}^{3} \vn \partial_i \varphi \wedge (A_i \vB),
\end{eqnarray*}	 
where we obtain the first and third term in (\ref{eq01}). \hfill$\blacksquare$\\

\noindent{\bf Acknowledgments:} J. \textsc{He} is supported by the program \emph{Sophie Germain} of the \emph{Fondation Math\'ematique Jacques Hadamard}.

\quad\\
\begin{multicols}{4}
\begin{minipage}[r]{50mm}
Diego \textsc{Chamorro}\\[3mm]
{\footnotesize LaMME\\UMR 8071\\UEVE\\Evry - France}
\end{minipage}
\begin{minipage}[r]{50mm}
Fernando \textsc{Cortez}\\[3mm]
{\footnotesize EPN\\ Quito - Ecuador}
\end{minipage}
\begin{minipage}[l]{50mm}
Jiao \textsc{He}\\[3mm]
{\footnotesize LaMME\\UMR 8071\\UEVE\\Evry - France}
\end{minipage}
\begin{minipage}[l]{50mm}
Oscar \textsc{Jarr\'in}\\[3mm]
{\footnotesize DIDE\\UTA\\Ambato - Ecuador}
\end{minipage}
\end{multicols}


\begin{thebibliography}{2}
\bibitem{Adams} 
D. R. \textsc{Adams} \& J. \textsc{Xiao}. \emph{Morrey spaces in harmonic analysis}. Ark. Mat. Volume 50, Number 2, 201-230 (2012).
\bibitem{ChenMiao}
Q. \textsc{Chen}, C. \textsc{Miao} \& Z. \textsc{Zhang}. \emph{On the Regularity Criterion of Weak Solution for the 3D Viscous Magneto-Hydrodynamics Equations}. Commun. Math. Phys. 284, 919–930 (2008).
\bibitem{CKN} 
L. \textsc{Caffarelli}, R. \textsc{Kohn} \& L. \textsc{Nirenberg}. \emph{Partial regularity of suitable weak solutions of the Navier--Stokes equations}. Comm. Pure Appl. Math., 35:771-831 (1982).
\bibitem{CML} D. \textsc{Chamorro}, K. \textsc{Mayoufi} \& P.-G. \textsc{Lemari\'e-Rieusset}. \emph{The role of the pressure in the partial regularity theory for weak solutions of the Navier-Stokes equations}. Archive for Rational Mechanics and Analysis, 228(1), 237-277. (2018)
\bibitem{Folland}
G. \textsc{Folland} \& E. \textsc{Stein}. \emph{Hardy spaces on homogeneous groups}. Princeton University Press, (1982). 
\bibitem{Chen}
Z.M. \textsc{Chen} \& W.R. \textsc{Price}. \emph{Morrey space techniques applied to the interior regularity problem of the Navier--Stokes equations}. Non-linearity; 14: 1453--1472 (2001).
\bibitem{Jia}
X. \textsc{Jia} \& Y. \textsc{Zhou}. \emph{Ladyzhenskaya–Prodi–Serrin type regularity criteria for the 3D incompressible MHD equations in terms of $3\times 3$ mixture matrices}. Nonlinearity, Volume 28, Number 9, (2015).
\bibitem{Kukavica} 
I. \textsc{Kukavica}. \emph{On partial regularity for the Navier--Stokes equations}. Discrete and continuous dynamical systems, 21:717-728 (2008).
\bibitem{Larios}
A. \textsc{Larios} \& Y. \textsc{Pei}. \emph{On the local well-posedness and a Prodi–Serrin-type regularity criterion of the three-dimensional MHD-Boussinesq system without thermal diffusion}. J. Differential Equations 263, 1419–1450, (2017).
\bibitem{OLeary}
M. \textsc{O'Leary}. \emph{Conditions for the local boundedness of solutions of the Navier--Stokes system in three dimensions}. Comm. Partial Differential Equations, 28:617-636 (2003).
\bibitem{PGLR2} P.G. \textsc{Lemarié-Rieusset}. \emph{Recent developments in the Navier-Stokes problem}. Chapman \& Hall/CRC, (2002).
\bibitem{Saka}
K. \textsc{Saka}. \emph{Besov Spaces and Sobolev spaces on a nilpotent Lie group}. Thoku. Math. Journ.  Vol. 31, p. 383-437 (1979).
\bibitem{Serrin1}
J. \textsc{Serrin}. \emph{On the interior regularity of weak solutions of the Navier--Stokes equations}. Arch. Rat. Mech. Anal., 9:187-195 (1962).
\bibitem{Struwe}
M. \textsc{Struwe}. \emph{On partial regularity results for the Navier--Stokes equations}. Comm. Pure Appl. Math., 41:437-458 (1988).
\bibitem{Taka}
S. \textsc{Takahashi}. \emph{On interior regularity criteria for weak solutions of the Navier--Stokes equations}. Manuscripta Math., 69:237-254 (1990).
\end{thebibliography}
\end{document}